\documentclass[11pt]{article}

\usepackage{setspace}
\usepackage{geometry}
\usepackage{amsmath}
\usepackage{graphicx}
\usepackage[colorlinks=true, allcolors=blue]{hyperref}
\usepackage{graphicx}
\usepackage{natbib}
\usepackage{amsmath} 
\usepackage{amssymb} 
\usepackage{placeins}

\usepackage{enumitem}
\usepackage{booktabs,threeparttable} 
\usepackage{appendix}

\usepackage{amsthm}
\theoremstyle{plain}

\theoremstyle{definition}

\theoremstyle{remark}

\usepackage{multirow}  

\newtheorem*{remark}{Remark}
\usepackage{graphicx}
\usepackage{multirow}
\usepackage{booktabs}

\usepackage[labelfont=bf,justification=raggedright,singlelinecheck=false]{caption}
\usepackage{xcolor}

\usepackage{float}  

\usepackage{authblk}

\makeatletter
\def\@fnsymbol#1{\ifcase#1\or *\or **\or ***\or §\or ¶\or **\or ††\or ‡‡\else\@ctrerr\fi}
\makeatother

\title{Specification tests for normal/gamma and stable/gamma  stochastic frontier models based on empirical transforms}

\author[1]{Christos K. Papadimitriou}
\author[1,2]{Simos G. Meintanis \thanks{Corresponding author: Simos G. Meintanis (on sabbatical leave from the University of Athens), email: simosmei@econ.uoa.gr } } 
\author[3]{Bernardo B. Andrade\thanks{On leave} }
\author[  ]{Mike G. Tsionas \thanks{Deceased}}

\affil[1]{Department of Economics, National and Kapodistrian University of Athens, Athens, Greece}
\affil[2]{Pure and Applied Analytics, North--WestUniversity, Potchefstroom, South Africa}
\affil[3]{Department of Statistics, University of Brasilia, Brasilia, Brasil}
\affil[  ]{ }

\date{}

\begin{document}
\maketitle

\begin{abstract}

Goodness--of--fit tests for the distribution of the composed error term in a Stochastic Frontier Model (SFM) are suggested.  The focus is on the case of a normal/gamma SFM and the heavy--tailed stable/gamma SFM. In the first case the moment generating function is used as tool while in the latter case the characteristic function of the error term is employed. In both cases our test statistics are formulated as weighted integrals of properly standardized data.  The new normal/gamma test is consistent, and is shown to have an intrinsic relation to moment--based tests. The finite--sample behavior of  resampling versions of both tests is investigated by Monte Carlo simulation, while  several real--data applications are also included. 
\newline 


\end{abstract}

\section{Introduction}
\label{intro}

Let $Y$ denote the maximum log--output obtainable from a (row) vector of log--outputs $X$, and consider the stochastic frontier production model (SFM),
\begin{equation} \label{SFM}
Y= X\beta+\varepsilon,  
\end{equation}
where $\beta$ is an unknown (column) vector of parameters, and $\varepsilon$ denotes the  error term. As it is standard practice in the context of SFMs, the error term is decomposed  as $\varepsilon=v-u$, where the random component $v$ is intended to capture purely random effects, while $u\geq0$ captures the technical inefficiency of the underlying firm. In this connection, the exponential and half/normal distribution are the most commonly used laws for the one-sided error component $u$. However more complicated distributions, such as the Gamma distribution have also been studied. 

Earlier work on goodness--of--fit tests for certain aspects of the SFM include Schmidt and Lin \citeyearpar{SL84},  Coelli \citeyearpar{Co95}, Lee \citeyearpar{Lee83}, Kopp and Mullahy \citeyearpar{KM90}, and Bera and Mallick \citeyearpar{BeMa02}. The work herein however is closer to Wang et al. \citeyearpar{WAS11} because as in that paper we also suggest omnibus tests, i.e. tests that have non--trivial power against arbitrary deviations from the null hypothesis, and use resampling in order to capture the extra variability introduced at the estimation of parameters step preceding the test. We also share common tools with Chen and Wang \citeyearpar{CW12} who use the characteristic function  for testing distributional specifications in the context of SFMs. 
Herein we follow the approach of Meintanis and Papadimitriou \citeyearpar{MP22} in suggesting specification tests for SFMs based on distributional transforms such as the moment generating function (MGF) and the characteristic function (CF). In this connection we introduce a two--level generalization  of the methods of  Meintanis and Papadimitriou  \citeyearpar{MP22}, first from the normal/exponential SFM to the normal/gamma SFM using the MGF, and then to the even more general stable/gamma SFM. Because the latter model is heavy--tailed we employ the CF rather than the MGF, both in the estimation step as well as in formulating our tests.  We finally mention a significant most recent contribution by Cheng et al. \citeyearpar{CWXZ} that propose versions of the classical distribution--function based Kolmogorov--Smirnov and Cram\'er--von Mises  goodness--of--fit tests for the law of the inefficiency component $u$. Their tests are more general than ours in that they do not assume a specific parametric form for the production function and that they are not  tailored, but apply to any distribution under test and allow for heteroskedasticity. On the other hand the tests of  Cheng et al. \citeyearpar{CWXZ}  can not be applied to our setting because they assume that all other parameters besides scale are fixed, while our tests are fully composite in the sense that all distributional parameters are estimated from the data at hand. We should also mention the test of Papadopoulos and Parmeter \citeyearpar{PP21} that depends on a combination of sample skewness and sample kurtosis. This test  is confined to identification of the one--sided error term and assumes the existence of higher moments of the error distribution, and therefore is not omnibus and can not be applied to the case of the stable/gamma SFM. Nevertheless it is quite simple to implement and moreover was shown to compete well, and in fact outperform, one of the non--omnibus  tests suggested by  Wang et al. \citeyearpar{WAS11}.

In our methods we employ as tools the MGF and the CF which are much more easier to apply in the situations considered herein, namely for a normal/gamma SFM as well as for the stable/gamma SFM because the corresponding distribution functions are quite complicated. A further advantage of the methods suggested herein is that they allow one to distinguish between the former model involving the normal distribution and a purely heavy--tailed stable/gamma SFM. 

The rest of the paper runs as follows. In Section \ref{sec_2} we introduce the test statistic for the normal/gamma SFM, show its consistency and discuss the role of the weight function. Simulation results for a resampling version of the test for the normal/gamma null hypothesis are presented in Section \ref{sec_4}. In Section \ref{sec_5} we consider the more general problem of testing the stable/gamma SFM. Simulations for this test are also included. Section \ref{sec_6} contains application of the proposed  tests  to several data sets from the literature.  Finally, conclusions and outlook are presented in Section \ref{sec_7}. Some technical material is deferred to \ref{ap1}--\ref{ap5}, while additional real--data applications are included in an on--line {\color{blue}{Supplement}}.

\section{Test for the normal/gamma SFM}\label{sec_2}

\subsection{Test statistic}\label{sec_2.1}

Consider the SFM in equation \eqref{SFM} where $\varepsilon=v-u$, where $u$ and $v$ are distributed independently of each other.  In this context and on the basis of data $(X_j,Y_j),  \ j=1,...,n$, we wish to test the null hypothesis
\begin{equation} \label{null}
\begin{split}
{\cal{H}}_{0}: \; & \textrm{Model (\ref{SFM}) holds true with }  
v \sim  {\cal{N}}(0,\sigma_v^2),  \ u\sim {\cal{G}}(p,c), 
 \textrm{ for some } \sigma_v^2, p,c >0,
\end{split}
\end{equation}
where ${\cal{N}}(\mu,\sigma^2)$, denotes a normal distribution with  mean $\mu$ and variance $\sigma^2$, and ${\cal{G}}(p,c)$ denotes a gamma distribution with shape parameter $p$ and scale parameter $c$. 

Recall in this connection that the moment generating function (MGF) of a generic random variable, say $Z$, is defined by $M_Z (t)=\mathbb E({e^{tZ}}), \: t\in \mathbb R$, while if $Z\geq 0$, then  its Laplace transform is given by $L_Z(t)=\mathbb E(e^{-t Z}), \ t>0$. For later use we also mention that if $Z\sim {\cal{N}}(\mu,\sigma^2)$ with density $(\sigma \sqrt{2\pi})^{-1}e^{-(1/2)((z-\mu)/\sigma)^2}$ then its MGF is given by  $M_Z (t)=e^{t\mu+ \frac{1}{2} \sigma^2 t^2}$, while if $Z \sim {\cal{G}}(p,c)$ with density $z^{p-1}/(c^p {\rm{\Gamma}}(p))e^{-z/c}, \ z>0$, then  $L_Z (t)= (1 + c t)^{-p}$, where ${\rm{\Gamma}}(p)=\int_{0}^{\infty} z^{p-1}e^{-z} {\rm{d}}z$ is the gamma function.

Based on the independence of $u$ and $v$,  it readily follows that under the null hypothesis ${\cal{H}}_0$, the MGF of the composed error term $\varepsilon$ may be computed as 
\begin{equation} \label{mgf15}
M_{\varepsilon} (t)=M_v(t)L_u(t)=\frac{e^{ \frac{1}{2}\sigma_v^2 t^2}}{(1+c t)^p}, \ t>0,
\end{equation}
and thus satisfies the differential equation 
\begin{equation} \label{mgf14}
(1+c t)M'_{\varepsilon}(t)+[c p-\sigma^2_v t(1+ct)]M_{\varepsilon}(t)=0 
\end{equation},
see also  Henze et al. \citeyearpar{EHM12}. 

Consequently  the MGF of the standardized composed error term $\widetilde \varepsilon=\varepsilon/c$ is given by 
\begin{equation} \label{mgf11}
M_{\widetilde\varepsilon}(t)=\frac{e^{ \frac{1}{2} \lambda t^2}}{(1+t)^p}, 
\end{equation}
and thus satisfies the differential equation 
\begin{equation} \label{mgf12}
(1+t)M'_{\widetilde\varepsilon}(t)+[p-\lambda t(1+t)]M_{\widetilde\varepsilon}(t)=0, \ t>0, 
\end{equation}
 where  $\lambda=\frac{\sigma^2_v}{c^2}$. 

Now let $\widehat{\varepsilon}_j=Y_j-\widehat\beta X_j$, be the residuals of the SFM  and write $r_j=\widehat{\varepsilon}_j/\widehat{c}, \ j=1,\ldots,n$, for the respective standardized residuals, involving suitable parameter estimates ($\widehat{\beta},\widehat{\sigma}_{v}^{2},\widehat c,\widehat p$), obtained under the null hypothesis ${\cal{H}}_0$. Then the left-- hand side of equation (\ref{mgf12}) may be estimated by 
\begin{equation} \label{mgf13}
D_n(t):=(1+t)M'_n(t)+[\widehat p-\widehat \lambda t(1+t)]M_n(t),  \ t>0 , 
\end{equation}
where 
\begin{equation} \label{emgf}
M_n(t)=\frac{1}{n} \sum_{j=1}^{n}e^{tr_j}, \ M'_n(t)=\frac{1}{n} \sum_{j=1}^{n}r_j e^{t r_j},
\end{equation}
stand for the empirical MGF based on $r_j$ and its derivative, respectively.
In view of equation \eqref{mgf12}, we should have for large $n$, $D_n(t)\approx 0$ under the null hypothesis, identically in $t$. Thus it is reasonable to reject the null hypothesis ${\cal{H}}_0$ of a normal/gamma SFM for large values of the test statistic
\begin{equation} \label{ts}
T_{n,w}= n \int_{0}^{\infty} D^2_{n}(t) \; w(t) \; {\rm{d}}t,
\end{equation}
where $w(t)>0$ is an integrable weight function.

Due to the  uniqueness of the solution of the differential equation figuring in \eqref{mgf14} and the uniqueness of MGFs, it follows that the test statistic $T_{n,w}$ defined by equation \eqref{ts} satisfies \[
T_{n,w} \rightarrow \infty, \ \mbox{a.s. as} \  n\to\infty, 
\]
under alternatives, and thus the test which rejects the null hypothesis ${\cal{H}}_0$ for large values of $T_{n,w}$ has asymptotic power one against arbitrary deviations from ${\cal{H}}_0$.  The argument leading to this property may be found in \ref{ap1}.

The choice $w(t)=e^{-\gamma t^2}, \gamma>0$, is particularly appealing from the computational point of view as it leads after straightforward algebra to the explicit formula

\begin{equation} \label{statistic}
\begin{split}
T_{n,\gamma} & = 
\frac{1}{n} \sum_{j,k=1}^{n}
\Bigg\{
\frac{r_jr_k+r_j(\widehat{p}-\widehat{\lambda})-\widehat{p}\widehat{\lambda}}{\gamma} +\frac{r_jr_k+\widehat{\lambda}\big(\widehat{\lambda}-4r_j-2\widehat{p}\big)}{4\gamma^2} (r_j+r_k) \\
& \hspace{5em} +\frac{\widehat{\lambda}(\widehat{\lambda}-r_j)}{4\gamma^3} \big[(r_j+r_k)^2+4\gamma\big] +\frac{\widehat{\lambda}^2}{16\gamma^4}(r_j+r_k)\big[(r_j+r_k)^2+10\gamma\big]   \\
& \hspace{5em}  +\frac{\sqrt{\pi}}{\sqrt{\gamma}}e^{\frac{(r_j+r_k)^2}{4\gamma}}
\;{\rm{\Phi}}\Big(\frac{r_j+r_k}{\sqrt{2\gamma}}\Big) 
\Bigg[ (r_j r_k+2\widehat{p}r_j+\widehat{p}^2)   +\frac{r_jr_k + r_j(\widehat{p}-\widehat{\lambda})-\widehat{p}\widehat{\lambda}}{\gamma}(r_j+r_k)  \\
& \hspace{5em} +\frac{r_jr_k+\widehat{\lambda}\big(\widehat{\lambda}-4r_j-2\widehat{p}\big)}{4\gamma^2} \big[(r_j+r_k)^2+2\gamma\big]  +\frac{\widehat{\lambda}(\widehat{\lambda}-r_j)}{4\gamma^3}(r_j+r_k)\big[(r_j+r_k)^2+6\gamma\big]  \\ 
& \hspace{5em} +\frac{\widehat{\lambda}^2}{16\gamma^4}\big[12\gamma^2+12\gamma(r_j+r_k)^2+(r_j+r_k)^4\big]
\Bigg]
\Bigg\}, 
\end{split}
\end{equation}
where ${\rm{\Phi}}(\cdot)$ denotes the standard normal distribution function. Notice that we write $T_{n,\gamma}$ for the test statistic in  \eqref{ts} corresponding to the weight function  $w(t)=e^{-\gamma t^2}$,  and thus make the dependence of the test statistic on $\gamma>0$ explicit.  Further details on the above computational formula for $T_{n,\gamma}$ are provided in   \ref{ap2}. 

\subsection{The role of the weight function}\label{sec_2.2}
Since the parameters of the model are unknown we need to efficiently estimate them from the data $(X_j,Y_j), \ j=1,\ldots,n$. Maximum likelihood estimators appear to be the most common, however drawbacks emerge when a closed form for the density of the random variable is not available as in the present case; see Tsionas \citeyearpar{Tsionas2012} for alternative estimation methods for SFMs whose likelihood function is not available in closed form. The corrected  least squares (COLS) estimation is the method employed herein. This method is further specified in \ref{ap3}.  

At this point, we will illustrate the role that the weight function $e^{-\gamma t^2}$ plays in the test statistic $T_{n,\gamma}$.
To this end, we compute the first five moments of the composed  normal/gamma error $\varepsilon$ in equation \eqref{SFM} as  
\begin{equation} \label{moments} 
\begin{split}
& \mu_1 = -c p, \mu_2=\sigma^2_v+c^2p(p+1),  \mu_3=-c^3 p(p+1)(p+2)-3cp\sigma^2_v, \\ 
& \mu_4 =3 \sigma^4_v+c^4 p(p+1)(p+2)(p+3)+6\sigma^2_v c^2 p(p+1), \\
& \mu_5=-15 \lambda^2 p-10 \lambda p (p+1)(p+2)-p(p+1)(p+2)(p+3)(p+4),
\end{split}
\end{equation}
where $\mu_k=\mathbb E(\varepsilon^k)$.

From equation \eqref{moments} we readily obtain that under the null hypothesis ${\cal{H}}_0$ of a normal/gamma SFM, the moments of the  standardized error $\widetilde \mu_k=\mathbb E(\widetilde \varepsilon^k)$ satisfy the  moment equations 
\begin{equation} \label{eq_moments}  {\cal{M}}_k=0, \ k=1,2,3,4,5,  \end{equation} where
\begin{equation} \label{moments_stand} 
\begin{split}
& {\cal{M}}_1= \widetilde \mu_1+p, \; {\cal{M}}_2=\widetilde \mu_2+(p+1)\widetilde\mu_1-\lambda, \; {\cal{M}}_3=\widetilde\mu_3+(p+2)\widetilde\mu_2-2\lambda(\widetilde \mu_1+1), \\
& {\cal{M}}_4=\widetilde \mu_4+(p+3)\widetilde \mu_3-3\lambda (\widetilde \mu_2+2\widetilde \mu_1),
\; {\cal{M}}_5=\widetilde \mu_5+(p+4)\widetilde \mu_4-4\lambda (\widetilde \mu_3+3\widetilde \mu_2).
\end{split}
\end{equation}

Now consider the expansion $e^z=1+z+z^2/2+z^3/6+z^4/24+{\rm{o}}(z^4), \ z\to 0$, of the exponential function figuring in \eqref{emgf}, and then starting from  \eqref{mgf13} and after some long but straightforward algebra we obtain
\begin{equation} \label{Dn} 
D_n(t)={\cal{M}}_{1,n}+t {\cal{M}}_{2,n}+\frac{t^2}{2}{\cal{M}}_{3,n} + \frac{t^3}{6}{\cal{M}}_{4,n}+\frac{t^4}{24}{\cal{M}}_{5,n}+{\rm{o}}(t^4), 
\end{equation}
$t\to 0$, where ${\cal{M}}_{k,n}$ are the empirical moment equations resulting from \eqref{moments_stand} if we replace the population moments $\widetilde \mu_k$ by the sample moments $\widetilde m_k=n^{-1} \sum_{j=1}^n r^k_j$ of the standardized residuals, $k=1,2,3,4,5$, and the parameters $(\lambda,p,c)$ by estimators $(\widehat \lambda, \widehat p, \widehat c)$.

A little reflection on equation \eqref{Dn} illustrates the role of the weight function and of the tuning parameter $\gamma$. Specifically the weight function provides the particular functional form of the weights by means of which each empirical moment equation figuring in \eqref{Dn} enters the test statistic $T_{n,\gamma}$. Consequently the expansion in \eqref{Dn} underlines the fact that a higher value of $\gamma$ results in a reduced impact of higher moments and vice versa, smaller $\gamma$ allows these higher moments to have a more significant impact on $T_{n,\gamma}$. In this connection caution should be exercised when choosing $\gamma$: A value that is too large allows only lower moments to enter the test statistic which in turn will adversely affect the power of $T_{n,\gamma}$, while on the other hand, a $\gamma$--value which is too small results in a test statistic prone to numerical error due to divergence of the exponential terms involved. Note in this connection that for $\gamma=0$, the value of  $T_{n,\gamma}$ diverges. On the other hand, as it is shown in \ref{ap3},  the collection of test statistics $\{T_{n,\gamma}, \ \gamma>0\}$ is closed at the upper boundary, i.e. the limit $\lim_{\gamma\to\infty} T_{n,\gamma}$ exists and is finite.


\section{Simulations}\label{sec_4} 

In this section we present the results of Monte Carlo study for the new test statistic given by equation \eqref{statistic}\footnote{For simulations we used Matlab software R2020b version.}. We consider the normal/gamma SFM, whereby $\varepsilon\sim {\cal{NG}}(\sigma^2_v,p,c)$ denotes a normal/gamma distribution with the indicated parameters. 

We consider the simple location SFM in equation (\ref{SFM}), $Y_j=\beta + \varepsilon_j$. 
The parameter vector $\Theta:=(\beta,\sigma^2_v,p,c)$ is estimated by the COLS method. The number of Monte Carlo replications is  $M=1,000$, with sample size $n=50, 100,200$, and nominal level of significance $5\%$. For reasons explained in Section \ref{sec_2.2}, we avoided to choose values of the tuning parameter $\gamma$ close to zero, choosing to implement our simulations for $\gamma=4$, $\gamma=6$ and $\gamma=8$.
 
Since the parameters of the model are considered unknown we employ a parametric bootstrap version of the tests in order to obtain critical values of the test statistic, say $T$ for simplicity. Specifically, we implement the so-called ``warp speed'' bootstrap procedure in order to approximate the bootstrap critical values in the Monte Carlo study (see Giacomini et al. \citeyearpar{GPW13} and Chang and Hall \citeyearpar{CH15}). This procedure is carried out by performing the following steps: 
\begin{itemize}
\item [(S1)] Draw a Monte Carlo sample $\{Y^{(m)}_j,X^{(m)}_j\}, \ j=1,...,n$,  compute the estimator--vector $\widehat{\Theta}^{(m)}$, where $\widehat{\Theta}^{(m)}=(\widehat{\beta}^{(m)}, \widehat{\sigma}^{2(m)}_v, \widehat{p}^{(m)},\widehat{c}^{(m)})$ 
\item[(S2)] On the basis of $\widehat{\Theta}^{(m)}$ calculate the residuals $\widehat{\varepsilon}^{(m)}_j$ and  the corresponding test statistic $T_m=T(\widehat{\varepsilon}^{(m)}_1,\ldots,\widehat{\varepsilon}^{(m)}_n)$
\item[(S3)] Generate {\it{i.i.d.}} bootstrap errors $\varepsilon_j^{(m)},\ j=1,\ldots,n$,  where $\varepsilon_j^{(m)}=v_j^{(m)}-u_j^{(m)}$, with $v_j^{(m)} \sim {\cal{N}}(0,\widehat{\sigma}_v^2)$ and $u_j^{(m)} \sim {\rm{Gamma}} (\widehat{p},\widehat{c})$, and independent.
\item[(S4)] Define the bootstrap observations $Y_j^{(m)}= \widehat{\beta}^{(m)}X_j+\varepsilon_j^{(m)}$, $j=1,...,n$.  
\item[(S5)] Based on $\{Y_j^{(m)},X_j\}$  compute the bootstrap estimator $\widehat{\Theta}_b^{(m)}=(\widehat{\beta}_b^{(m)}, \widehat{\sigma}_{b,v}^{2 (m)}, \widehat{p}_b^{(m)},\widehat{c}_b^{(m)})$, and the corresponding bootstrap residuals, say, $\widehat{\epsilon}_{j}^{(m)}$, $j=1,...,n$. 
\item[(S6)] Compute the test statistic $\widehat T_m:=T(\widehat{\epsilon}^{(m)}_{1},\ldots,\widehat{\epsilon}^{(m)}_{n})$ 
\item[(S7)] Repeat steps (S1)--(S6), for  $m=1,...,M$, and collect  test--statistic values $T_m$ and bootstrap statistic values $\widehat T_m, \ m=1,...,M$. 
\item[(S8)] Set the critical point, say $C$, equal to $\widehat T_{(M-0.05 M)}$, where $\widehat T_{(m)}$, $m=1,\ldots,M$, denote the order statistics corresponding to $\widehat{T}_{m}$.
\item[(S9)] Obtain the empirical rejection rate $M^{-1} \sum_{m=1}^M {\bf{1}}_{T_{m}>C}$.
\end{itemize}

In \autoref{clslevel}  the size results are presented for the test statistic $T_{n,\gamma}$, corresponding to the ${\cal{NG}}(1,p,1)$ null hypothesis for $p={\color{blue}0.25, 0.5}, 1,2,3$. The power of the test is also obtained against an alternative mixture hypothesis $0.7{\cal{NG}}(1, 1, 1)+0.3{\cal{NG}}(1, p,1)$,   whereby an observation is drawn from a ${\cal{NG}}(1, 1, 1)$ distribution (resp. ${\cal{NG}}(1, p, 1)$ distribution) with probability 0.7 (resp. 0.3), for $p={\color{blue}0.25, 0.4}, 0.5, 2, 3$; see \autoref{clspower}.

\begin{table}[htbp]
\caption{Size of the test for the Normal/Gamma null hypothesis ${\cal{NG}}(1,p,1)$ with sample size $n$. Nominal level of significance $5\%$.}
\label{clslevel}\centering
\setlength{\tabcolsep}{20pt} 
\renewcommand{\arraystretch}{0.58}
\small
{%
\begin{tabular}{@{}cc  ccc@{}}   
& $n$ & $\gamma$=4 & $\gamma$=6 &$\gamma$=8 \\
\midrule 
$p$=0.25 & 50 & 1.5 & 1.9 & 2.0 \\
    & 100  & 3.8 & 3.4 & 4.3 \\
    & 200  & 2.5 & 3.2 & 3.7 \\
    & 400  & 3.0 & 2.8 & 3.3 \\
\midrule
$p$=0.5 & 50 & 3.3 & 3.2 & 3.9 \\
    & 100  & 6.0 & 5.0 & 5.0 \\
    & 200  & 3.4 & 2.9 & 3.0 \\
    & 400  & 2.9 & 3.2 & 3.8 \\
\midrule
$p$=1.0 & 50 & 5.5 & 6.3 & 6.1 \\
    & 100  & 5.1 & 4.6 & 4.4 \\
    & 200  & 5.2 & 4.7 & 5.7 \\
    & 400  & 3.4 & 2.5 & 2.7 \\
\midrule
$p$=2.0 & 50 & 3.9 & 3.8 & 3.4	 \\
    & 100	& 5.7 & 6.2 & 6.2	 \\
    & 200	& 6.0 & 6.7 & 7.1 \\
    & 400  & 4.2 & 4.7 & 4.8 \\
\midrule
$p$=3.0 & 50 & 5.5 & 6.3 & 6.1 \\
    & 100  & 5.1 & 4.6 & 4.4 \\
    & 200  & 5.2 & 4.7 & 5.7 \\
     & 400  & 4.8 & 4.6 & 5.1 \\
\end{tabular}}
\end{table}

\begin{table}[htbp]
\caption{Power of the test against mixture alternatives of the form $0.7 {\cal{NG}} (1 ,1 ,1) + 0.3 {\cal{NG}} (1 ,p,1)$ with sample size $n$. }
\label{clspower}\centering
\setlength{\tabcolsep}{20pt} 
\renewcommand{\arraystretch}{0.58}
\small
{%
\begin{tabular}{@{}cc ccc@{}}   
& $n$ & $\gamma$=4 & $\gamma$=6 &$\gamma$=8 \\
\midrule 
$p=0.25$ 
& 50	& 19.0 & 24.5 & 26.5 \\
& 100	& 28.3 & 33.1 & 34.5 \\
& 200	& 30.6 & 38.3 & 40.5 \\ 
\midrule 
$p=0.40$ 
& 50	& 9.6 & 14.1 & 14.5 \\
& 100	& 29.3 & 31.4 & 32.7 \\
& 200	& 31.1 & 29.5 & 32.9 \\ 
\midrule 
$p=0.5$ 
& 50	& 12.2 & 14.3 & 14.9 \\
& 100	& 17.9 & 19.8 & 20.0 \\
& 200	& 38.1 & 46.8 & 49.2 \\ 
\midrule 
$p=2.0$     
& 50	& 33.0 & 35.1 & 35.3  \\
& 100	& 50.5 & 52.5 & 52.9  \\
& 200	& 62.3 & 65.4 & 68.1  \\
\midrule 
$p=3.0$     
& 50	& 55.3 & 54.9 & 55.3 \\
& 100	& 61.2 & 62.1 & 62.7 \\
& 200	& 74.4 & 75.7 & 76.5 \\ 

\end{tabular}}
\end{table}

From \autoref{clslevel} we observe that the empirical size of the test $T_{n,\gamma}$ varies with the value of the tuning parameter $\gamma$, but it is overall,  with a few exceptions, close to the $5\%$ nominal size of the test. Moreover, looking at \autoref{clspower} we observe that the power of the test progressively increases along with the sample size $n$. Also the value of the tuning parameter $\gamma$ is seen to exert a certain effect on the power of the test, which is however mostly mild.  

We have also examined whether the test is sensitive to violations of the pure random error distribution. In this connection we applied the test on a Student--t/gamma distribution, whereby the normal distribution in the normal/gamma SFM is replaced by a  Student--t distribution. In   \autoref{clspowertStG3} the power results are shown and correspond to $\nu$ degrees of freedom and $p=3$ for the gamma shape parameter. These results suggest that the test may not be particularly sensitive to violations of the normal/gamma null hypothesis  towards alternative  laws for the pure random error $v$.

\begin{table}[htbp]
\caption{Power of the test against a  Student--t/Gamma alternative with $\nu$ degrees of freedom and $p=3$ with sample size $n$.}
\label{clspowertStG3}\centering
\setlength{\tabcolsep}{20pt} 
\renewcommand{\arraystretch}{0.70}
\small
{%
\begin{tabular}{@{}cc  ccc ccc@{}}   
& $n$ & $\gamma$=0.5 & $\gamma$=1 &$\gamma$=2& $\gamma$=4 & $\gamma$=6 &$\gamma$=8 \\
\midrule
$\nu$=5 & 50	&11.6&11.0&10.3& 9.7 & 9.5 & 9.3 \\
    & 100	&11.5&10.1&11.7& 9.3 & 8.2 & 7.3 \\
    & 200	&12.5&10.1&12.4&11.6&9.5&8.5\\
\midrule
$\nu$=6 & 50 & 7.4 & 6.8 & 7.9 & 9.7 & 9.1 & 8.6 \\
    & 100  & 14.1 & 15.7 & 13.4 & 12.2 & 11.0 & 9.3 \\
    & 200  & 17.5 & 14.6 & 15.1& 13.1 & 10.7 & 8.3\\
\end{tabular}}
\end{table}

Finally, we point out that we have obtained similar results with maximum likelihood estimation (see next section) instead of COLS with, on average, a decrease of 0.5 percentage points in the empirical sizes compared to \autoref{clslevel}.

In this connection and as already implicit, the weight function $w(t)$ figuring in equation \eqref{ts} may in principle take arbitrary functional forms. A trivial condition that leads to test consistency  is that $w(t)$ should be positive, while a further requirement already encountered in the derivation of equation \eqref{statistic} is computational convenience, meaning that it should render the test statistic in a closed--formula free of numerical integration. The most popular weight function so far has been the zero--mean normal density, that is also adopted herein. There exist some technical works on how to properly choose the weight parameter $\gamma$, in the case of the zero--mean normal density as weight function. It should be kept in mind though that this choice is distribution-- and model--specific and thus these results only apply to the specific context under which they were derived.    

In order to see where the problem lies we note that the limit null distribution of $T_{n,w}$ is that of $\sum_{j=1}^\infty \omega_j {\cal{N}}^2_j$, where $({\cal{N}}_j, \ j\geq 1)$ are i.i.d. standard normal random variables. In turn $\omega_1\geq \omega_2\geq \cdots$, are the eigenvalues of a complicated integral equation depending on the weight function $w(\cdot)$, the distribution being tested as well as on the estimation of any parameters involved. Regarding power properties, the notion of Bahadur efficiency is often the preferred notion, but such efficiency requires calculation of the eigenvalues, and in particular the knowledge of $\omega_1(=\omega_1(\gamma)$ if $w(t)={\rm{e}}^{-\gamma t^2}$), the largest eigenvalue, which however is rarely available. The corresponding technical analysis is highly non--trivial and in general calculation of eigenvalues  remains an open problem. Nevertheless some results in this direction are available in \cite{ten1}, \cite{ten2}, \cite{mmo} and \cite{eh}. Here instead we adopt the pragmatic approach of observing how the finite--sample power of the test varies with the parameter $\gamma$,  thus making empirically--based choices. In this regard we refer to our Monte Carlo results which indicate that the power varies to some extend with $\gamma$, with the value $\gamma=8$ leading to somewhat more powerful test. Nevertheless our results are not thorough and more work is needed in this direction, always keeping in mind that the choice of $\gamma$ depends on the distribution being tested but also on the direction away from the null hypothesis, and that, according to the general results of  \cite{esc}, no universally optimal test exists, however elaborate the choice of the weight parameter $\gamma$ may be.


\section{Testing for the stable/gamma model}\label{sec_5}
\subsection{Test statistic}\label{sec_5.1}
We now extend our test from the normal/gamma to the more general case of testing the stable/gamma SFM whereby the two--sided error $v$ follows a symmetric stable distribution, rather than a normal distribution. Specifically we wish to test the null hypothesis
\begin{equation}
\begin{split}
{\Upsilon_{0}}: \; & \textrm{Model (\ref{SFM}) holds true with }  
v \sim  {\cal{S}}(\alpha,\kappa),  \ u\sim {\cal{G}}(p,c), 
 \textrm{ for some } 0<\alpha\leq 2, \textrm{and} \ \kappa, p,c >0,
\end{split}
\end{equation}
where ${\cal{S}}(\alpha,\kappa)$, denotes a symmetric stable distribution with  tail--index $\alpha$ and scale parameter $\kappa$.

In this connection recall that the CF of a generic random variable $Z$ is defined by $\varphi_Z (t)=\mathbb E({e^{{\rm{i}}tZ}}), \: t\in \mathbb R$, ${\rm{i}}=\sqrt{-1}$, and that the CF of a symmetric stable distribution is given by $e^{-\kappa^\alpha|t|^\alpha}$\footnote{When we write stable distribution, we actually mean the symmetric stable distribution}. Note that for $\alpha=2$, this distribution reduces to the normal distribution with variance equal to $2\kappa^2$, which is the only stable distribution with finite variance. Otherwise if $\alpha<2$, the stable distribution has finite moments only of order strictly less than $\alpha$. As a consequence the MGF of $\varepsilon$ does not exist for any $t\neq0$, and thus we are led to use the CF as our testing tool, given that the CF also determines uniquely the underlying law of error. Note that the mean of the distribution exist only if  $\alpha>1$, in which case it is equal to zero. Accordingly, and since stable laws have potential for application mostly, if not exclusively, in the finite--mean situation, hereafter we only consider the case of $\alpha>1$. For further properties of this distribution we refer to Paolella \citeyearpar{Paolella07}. 

In order to introduce our test we note that under the  null hypothesis of a stable/gamma SFM, the CF of the standardized composed error term $\widetilde \varepsilon=\varepsilon/c$ is given by  
\begin{equation*} \label{cf21}
\varphi_{\widetilde\varepsilon} (t)=\frac{e^{ -\lambda |t|^\alpha}}{(1+{\rm{i}}t)^p},  
\end{equation*}
and consequently it satisfies the differential equation 
\begin{equation*} \nonumber
\Delta(t):=(1+{\rm{i}}t)\varphi'_{\widetilde\varepsilon}(t)+\left[{\rm{i}}p+\alpha \lambda |t|^{\alpha-1}(1+{\rm{i}}t){\rm{sgn}}(t)\right]\varphi_{\widetilde\varepsilon}(t)=0,  
\end{equation*}
where $\lambda=\left(\kappa/c\right)^\alpha$, and ${\rm{sgn}}(t)=-1, 0, 1$ if $t<0$, $t=0$, and $t>0$, respectively; see also Meintanis \citeyearpar{Meintanis05}. 
Thus arguing as before we suggest the test statistic
\begin{equation} \label{test3}
{\cal{T}}_{n,w}= n \int_{-\infty}^{\infty} \left|\Delta_{n}(t)\right|^2 \; w(t) \; {\rm{d}}t,   
\end{equation}
where   \[ 
\Delta_{n}(t)=(1+{\rm{i}}t)\varphi'_{n}(t)+\left[{\rm{i}}\widehat p+\widehat \alpha \widehat \lambda |t|^{\widehat\alpha-1}(1+{\rm{i}}t){\rm{sgn}}(t)\right]\varphi_{n}(t),
\]
with 
\[ \varphi_n(t)=\frac{1}{n} \sum_{j=1}^n e^{{\rm{i}}t r_j}:=C_n(t)+{\rm{i}} S_n(t)\]
and 
\[ \varphi'_n(t)=\frac{1}{n} \sum_{j=1}^n {\rm{i}} \: r_j e^{{\rm{i}}t r_j}
 :=C'_n(t)+{\rm{i}} S'_n(t)\]
being the empirical CF and its derivative, respectively, where \[C_n(t)=\frac{1}{n}\sum_{j=1}^n \cos t r_j, \ S_n(t)=\frac{1}{n}\sum_{j=1}^n \sin t r_j,\]  denote the real and imaginary part, respectively, of $\varphi_n(t)$, and $C'_n(t)$ and $S'_n(t)$ the corresponding parts of $\varphi_n'(t)$. Also we write     
$|z|=\sqrt{z^2_1+z^2_2}$, for the modulus of a complex number $z=z_1+{\rm{i}}z_2$.  The consistency of the test that rejects the null hypothesis of a stable/gamma SFM for large values of the test statistic ${\cal{T}}_{n,w}$ figuring in \eqref{test3}  may be shown by analogous arguments as in the case of the test in equation \eqref{ts}; see also Meintanis \citeyearpar{Meintanis05}. 

Following straightforward computations we obtain 

\begin{equation}  \label{d2}
\begin{split}
|\Delta_n(t)|^2 &= (1+t^2)|\varphi_n'(t)|^2  +\left[ (\widehat \alpha \widehat \lambda |t|^{\widehat\alpha-1})^2+ (\widehat\alpha \widehat\lambda |t|^{\widehat\alpha}+\widehat p)^2\right]|\varphi_n(t)|^2 + 2\widehat p \left(S'_n(t)C_n(t)-C'_n(t) S_n(t)\right) \\ 
&  \hspace{1em} +2\left[ \widehat \alpha \widehat \lambda |t|^{\widehat\alpha-1} {\rm{sgn}}(t)+t (\widehat p+\widehat \alpha\widehat\lambda |t|^{\widehat\alpha}) \right] (C_n'(t)C_n(t)+S_n'(t)S_n(t)) .
\end{split}
\end{equation}

Now if we consider $w(t)=e^{- \gamma t^2}$ as weight function, then the test statistic may be explicitly written by making use of the integrals
\[
I_{\nu,\gamma}(z)=\int_{0}^\infty t^{\nu} \cos(t z) e^{- \gamma t^2} {\rm{d}}t,  \ \nu >-2
\]
and
\[
J_{\nu,\gamma}(z)=\int_{0}^\infty t^{\nu} \sin(t z) e^{- \gamma t^2} {\rm{d}}t, \ \nu >-1. 
\]   

Specifically following a long but otherwise straightforward computation we have from \eqref{test3} and \eqref{d2},   
\begin{eqnarray} \label{test4} \nonumber
{\cal{T}}_{n,\gamma} & = & \frac{1}{n}\sum_{j,k=1}^{n}  \Big\{  2r_jr_k I_{2,\gamma}(r_j-r_k) + 2\big(r_jr_k + 2\widehat p r_j + \widehat p^2\big) I_{0,\gamma}(r_j-r_k) 
+ 4 \widehat p \widehat \alpha \widehat \lambda I_{\widehat \alpha ,\gamma}(r_j-r_k) \\ \nonumber
&+& 4 \widehat p r_k J_{1,\gamma}(r_j-r_k) 
 +  2\widehat \alpha^2\widehat \lambda^2\Big(I_{2\widehat\alpha,\gamma}(r_j-r_k)+I_{2(\widehat\alpha-1),\gamma}(r_j-r_k)\Big) \\ &+& 4\widehat\alpha\widehat\lambda r_k\Big(J_{\widehat\alpha+1,\gamma}(r_j-r_k)+J_{\widehat\alpha-1,\gamma}(r_j-r_k) \Big\} .
\end{eqnarray}

Further details on the above computational formula for  ${\cal{T}}_{n,\gamma}$ and for the integrals $I_{\nu,\gamma}(\cdot), J_{\nu,\gamma}(\cdot)$, are provided in  \ref{ap4} and  \ref{ap5}, respectively.

\subsection{Parameter estimation}\label{sec_5.2}
In a stochastic frontier model lacking a closed form for the likelihood,
we may estimate the parameters by means of the CF, following Tsionas \citeyearpar{Tsionas2012}
and Andrade and Souza \citeyearpar{Andrade2018}.
This can be done as described in  \ref{ap3}. Specifically we can approximate the density of the composed error, $f_\varepsilon$, by numerical inversion of its CF and proceed to likelihood maximization, 
which, in the context of the stable/gamma model, amounts to 
\begin{equation*} \label{estpar0}
\widehat \Psi :=  \underset{\Psi}{\mathrm{argmax}} \sum_{i} \ell_i(\Psi),
\end{equation*}
where $\ell_i(\Psi) := \log \hat{f}_{\varepsilon_i}(\cdot)$ 
(recall equation \ref{eq:dft}) and $\Psi:=(\beta,\kappa,\alpha, p,c)$.


\subsection{Monte Carlo study}\label{sec_5.3}

In this section we present the results of Monte Carlo study for the test statistic given by equation \eqref{test4}. We consider the stable/gamma SFM, whereby $\varepsilon\sim {\cal{SG}}(\kappa,\alpha, p,c)$ denotes a stable/gamma distribution with the indicated parameters. Again, we consider the simple location SFM in equation (\ref{SFM}), $Y_j=\beta + \varepsilon_j$. The parameters are estimated by maximum likelihood. We use the resampling of Section \ref{sec_4} with  $M=10,000$ Monte Carlo replications throughout. The sample sizes are $n$ = 200, 400, 500 for simulations under the nominal level of significance $5\%$ (Table \ref{cf_ml_level}). These sample sizes are larger than those used in Section \ref{sec_4} since stable distributions are known for requiring larger sample sizes for reliable parameter estimation. For power simulations we used $n=200$ and 500 (Table \ref{clspowerFFT}). The tuning parameter has been set to $\gamma=4$ and results are similar in the range $\gamma=4$ to $\gamma=8$ (not shown). Since our test seem to be rather conservative in estimating the nominal level, power results shown in Table \ref{clspowerFFT} have been corrected following  Lloyd \citeyearpar{Lloyd05} so that 
$$
\text{corrected power} = \Phi\left(\Phi^{-1}(pow) - 
\Phi^{-1}(\widehat{s}) + \Phi^{-1}(0.05) 
\right),
$$
where $pow$ is uncorrected (Monte Carlo) power and $\widehat{s}$ the (Monte Carlo) estimated size.

\begin{table}[htbp]
\caption{Size of the test for the Stable/Gamma null hypothesis ${\cal{SG}}(1,\alpha,1,1)$ with sample size $n$. }
\label{cf_ml_level}\centering
\setlength{\tabcolsep}{20pt} 
\renewcommand{\arraystretch}{0.70}
\small
\begin{tabular}{@{}cc cccc@{}}   
& $n$ & $\gamma$=2 & $\gamma$=4 & $\gamma$=6 &$\gamma$=8 \\
\midrule 
$\alpha$=1.8 & 200 &  1.4 & 3.9  & 3.2 & 2.1\\
    & 400  & 2.0 &  2.9 & 4.5 & 2.9\\
    & 500  & 3.3 & 4.9 & 4.6 & 3.5\\
\midrule
$\alpha$=1.9 & 200	& 3.8  & 1.1 & 3.0	& 1.0\\
    & 400	& 2.5 & 2.5 & 2.3&	2.7 \\
    & 500	& 1.1 & 3.3 & 2.4 & 1.6\\
\midrule
$\alpha$=1.95 & 200 &  3.6&  2.6&  2.5&  2.3\\
    & 400  & 1.5  & 3.1 & 3.4 & 3.7 \\
    & 500  & 1.9 & 3.4 & 3.3 & 1.7\\
\end{tabular}
\end{table}

\begin{table}[htbp]
\caption{Power (corrected) of the test with null hypothesis $\alpha = \alpha_0$ against different alternatives with sample sizes 200 and 500 and $\gamma=6$. }
\label{clspowerFFT}\centering
\setlength{\tabcolsep}{20pt} 
\renewcommand{\arraystretch}{0.75}
\small
{%
\begin{tabular}{@{}cc  cc@{}}   
alternative &  $n$ & $\alpha_0$=1.8 & $\alpha_0$=1.95 \\
\midrule 
$\alpha=1.5$ & 200 & 58 & 59  \\
&500 & 55 & 63\\
\midrule 
$\alpha=1.7$ &200 & 50 &  58 \\
& 500 & 51 & 59\\
\midrule 
$\alpha=1.8$ & 200 & $-$  & 52  \\
& 500 & $-$ & 56\\
\midrule 
$\alpha=1.9$ &200 & 52 &   48 \\
&500 &  45  &  50 \\
\midrule 
$\alpha=1.95$& 200 & 52  & $-$   \\
& 500 &  56  &  $-$ \\
 \midrule 
$t$(2)/Gamma& 200  & 53 &   53\\
& 500  & 54 &   56\\
\end{tabular}}
\end{table}

It should be noted that power simulations are trickier than in the normal/gamma case since computational time using maximum likelihood is much higher. In order to speed up simulations, we have fixed the value of the shape parameter at $\alpha=\alpha_0$, and simply tested against alternative values of $\alpha$ in $\{1.5,1.7,1.8,1.9,1.95\}$ under the null $\alpha_0=1.8$ and $\alpha_0=1.95$. We used sample sizes $n=200$ and $n=500$. 

Power values are quite high, ranging from  45\% and 63\%, but on the other hand sensitivity of the test with respect to sample size is rather limited. Specifically empirical rejection rates seem to reach a plateau already at $n=200$ and do not  further increase with the sample size, at least significantly so. Such behaviour was also observed in \autoref{clspowertStG3} and may be attributed to the fact that the distributions involved are heavy--tailed.

\begin{table}[htbp]

\end{table}


\section{Data and empirical results}\label{sec_6}

\subsection{Data}\label{sec_6.1}

In order to examine the behaviour of the proposed test statistics defined in  equations  \eqref{statistic} and  \eqref{test4} we use datasets from previous studies -- which are also employed by Horrace and Wang \citeyearpar{HW22} -- whereby we consider the normal/gamma and the stable/gamma  distribution as the null hypotheses of the composed error term. The datasets are in Christensen and Greene \citeyearpar{CG76} \footnote{The Christensen and Greene (1976) data are more recently analyzed by Tsionas (2002). These data are selected because they are available on the journal's website.} on US electric utility companies, in Tsionas \citeyearpar{Tsionas2006}  on US bank costs (considering a cost function), and in Polachek and Yoon \citeyearpar{PY96} on the earnings of workers (considering a production function). 
These data are publicly available in the \textit{Journal of Applied Econometrics Data Archive.}\footnote{\url{http://qed.econ.queensu.ca/jae/}}

In order to apply the test we consider the parametric bootstrap version of the test with $B=100$ bootstrap samples from the hypothesized normal/gamma (resp. stable/gamma) law with estimated parameters obtained  by COLS (resp. MLE). (We refer to 
Cornea--Madeira and Davidson \citeyearpar{CD15} for the validity of the parametric bootstrap in the case of heavy--tailed distributions such as the stable distributions). The nominal level of the test is set to $5\%$ and use $\gamma=1$ as value of the tuning parameter.\footnote{We also tried $\gamma=0.5$ and $\gamma=2$, leading to the same conclusions}

For consistency and comparison reasons we follow Tsionas \citeyearpar{Tsionas2002} and take the negatives of $Y$ and $X$ in the SFM of equation \eqref{SFM}.
\subsection{Empirical results}\label{sec_6.2}
First we turn our attention to a widely used dataset in stochastic frontier context, that of the cross-sectional case of 123 US electric utility companies in 1970 used by Christensen and Greene \citeyearpar{CG76} and Tsionas \citeyearpar{Tsionas2002}. We consider the following functional form for the cost function 
\begin{eqnarray*}\label{dataeq1}
-\log{\frac{C_j}{PF_{j}}}&=&-\beta_{0j}-\beta_{1j}\log{Q_j}-\beta_{2j}\log^2{Q_j} -\beta_{3j}\log{\frac{PL_{j}}{PF_{j}}}-\beta_{4j}\log{\frac{PK_{j}}{PF_{j}}}+v_j+u_j ,
\end{eqnarray*}
where $Q_j$ and $C_j$ are the output and the cost of the $j^{\rm{th}}$ firm respectively, and $PL,PK,PF$ denote the prices of the three production factors, namely labour, capital, and fuel. 

The parameter estimates, as well as the p-value of the stable/gamma test statistic are presented in \autoref{datapower1ML} using the ML estimation method described in  \ref{ap3}. The slope and intercept estimates are very close to those of Table I in Tsionas \citeyearpar{Tsionas2002}, while the p-value of the test statistic is clearly non--significant for a stable/gamma null hypothesis, meaning that we can not reject $\Upsilon_0$.  At the same time though, the estimate of the tail--index $\alpha$ borders the case of normality $\alpha=2$, and thus our findings  may be interpreted as being in agreement with those presented by Horrace and Wang \citeyearpar{HW22} implying a medium--tailed normal/gamma  SFM. This conclusion is further corroborated by an application of the normal/gamma test (with COLS estimation) on these data that leads to a significant increase of the resulting p--value; refer to the figures in the right panel of  \autoref{datapower1ML}.

\begin{table}
\centering
\caption{
Parameter estimates and p-value of the Stable/Gamma (resp. the Normal/Gamma)  test for the dataset in Christensen and Green (1976) with ML (resp. COLS) estimation.} 
\label{datapower1ML}. \centering
\setlength{\tabcolsep}{10pt} 
\renewcommand{\arraystretch}{1.00}
\small
{%
\begin{tabular}{  cc  cc  cc  cc  } 
 \multicolumn{4}{c}{Stable/Gamma} & \multicolumn{4}{c}{Normal/Gamma}  \\   
\toprule
$\widehat{\beta}_0$ & -6.790 & $\widehat{\kappa}$ & 0.0682 & \;\;\;\; $\widehat{\beta}_0$ & -7.205 &  & \\   
$\widehat{\beta}_1$ & 0.233 & $\widehat{\alpha}$ & 1.9941 & \;\;\;\; $\widehat{\beta}_1$ & 0.386 & $\widehat{\sigma}_v$ & 0.1409  \\
$\widehat{\beta}_2$ & 0.041 & $\widehat{p}$ & 0.3412 & \;\;\;\; $\widehat{\beta}_2$ & 0.032 & $\widehat{p}$ &   0.1700 \\ 
$\widehat{\beta}_3$ & 0.285 & $\widehat{c}$ & 0.2107 & \;\;\;\; $\widehat{\beta}_3$ & 0.246 & $\widehat{c}$ &   0.3987\\ 
$\widehat{\beta}_4$ & 0.121 & \multicolumn{2}{c}{{\bf{p--value}}=0.67} & \;\;\;\; $\widehat{\beta}_4$ & 0.079 & \multicolumn{2}{c}{{\bf{p--value}}=0.98} \\
\end{tabular}}
\end{table}

\begin{figure}[H]
\centering
\includegraphics[scale=0.8]{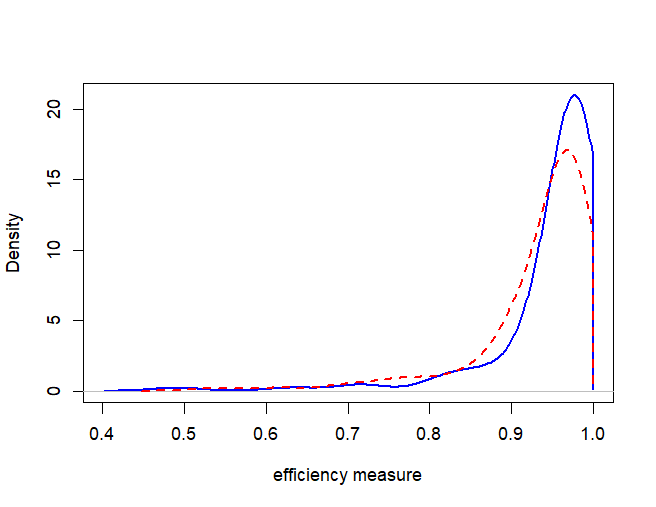}
\caption{Densities for efficiency measures, 
$\mathbb E(\exp(-u_i)|\varepsilon_i)$,
$i=1,\ldots,123$,
for two SFMs fitted to data from 
Christensen and Green (1976): normal/gamma (blue-full) and 
stable/gamma (red-dashed).}\label{fig:eff}
\end{figure}

 As detailed in  \ref{sec:eff-meas}, efficiency measures 
 can be obtained for the SFMs being considered in this paper. Efficiency 
 scores resulting from fitting the normal/gamma and the stable/gamma
 models to the 123 electric utility companies are depicted in Figure \ref{fig:eff}. In general, the normal/gamma model yielded slightly 
 more cost-efficient firms than the stable/gamma SFM.

We next revisit the panel data set in Polachek and Yoon \citeyearpar{PY96}, the \textit{Panel Study of Income Dynamics} (PSID), with the following specification of the production function
\begin{eqnarray*}\label{dataeq2}
\log W_{jt}&=&\beta_{0j}+\beta_{1j}EDU_{jt}+\beta_{2j}EX_{jt}+\beta_{3j}EX_{jt}^{2} +\beta_{4j}TEN_{jt}+\beta_{5j}TEN_{jt}^{2} + v_{jt}-u_{jt} ,
\end{eqnarray*}
where $W_{jt}, EDU_{jt}, EX_{jt}, TEN_{jt}$ represent the wage, education and experience in years and tenure in months, respectively, for $838$ workers.  
Corresponding results for the normal/gamma null hypothesis are presented in \autoref{datapower2} for the years 1970 and 1971. These results clearly lead to non--rejection of the null hypothesis ${\cal{H}}_0$ of a medium--tailed normal/gamma SFM.  This conclusion is somewhat inline with the corresponding results in Table ST.1 of the on-line Supplement where the p-value  leads to rejection of the null hypothesis $\Upsilon_0$ of a stable/gamma SFM for the year 1970, while for 1971, the estimated value of the stable shape parameter $\alpha$ borders that of normality ($\alpha=2$).

\begin{table}[H]
\caption{Parameter estimates and p-value of the Normal/Gamma test for the dataset in Polacheck and Yoon (1996) with COLS estimation}
\label{datapower2} \centering
\setlength{\tabcolsep}{10pt} 
\renewcommand{\arraystretch}{0.70}
\small 
{%
\begin{tabular}{l  cccc  l  cccc }  
\toprule
\multirow{6}{*}{1970} & $\widehat{\beta}_0$ & -1.2782 & 
& & \multirow{6}{*}{1971} & $\widehat{\beta}_0$ & -0.8086 &  &   \\
& $\widehat{\beta}_1$ & 0.0609 & $\widehat{\sigma}_v$ & 0.2654 & & $\widehat{\beta}_1$ & 0.067 & $\widehat{\sigma}_v$ & 0.1973  \\
& $\widehat{\beta}_2$ & 0.0603 & $\widehat{p}$ & 0.0897 & & $\widehat{\beta}_2$ & 0.0382 & $\widehat{p}$ & 0.5783 \\
& $\widehat{\beta}_3$ & -0.0012 & $\widehat{c}$ & 1.3690 & & $\widehat{\beta}_3$ & -0.0007 & $\widehat{c}$ & 0.4179 \\
& $\widehat{\beta}_4$ & 0.0069 & \multicolumn{2}{c}{{\bf{p--value}}=1.0} & & $\widehat{\beta}_4$ & -0.0067 & \multicolumn{2}{c}{{\bf{p--value}}=0.98} \\
& $\widehat{\beta}_5$ & -4.7E-5 & & & & $\widehat{\beta}_5$ & 5.2E-5 & &  \\ 
\bottomrule
\end{tabular}}
\end{table}

In the third real--data example, we employ the dataset of large\footnote{The commercial banks included in the Report of Condition and Income (Call Report) whose total assets exceeded one billion dollars for at least three (not necessarily consecutive) years.} US banks' costs used in Tsionas \citeyearpar{Tsionas2006}. The cost function is given as follows 
\begin{equation*} \label{dataeq3}
-\log C_{jt}= X_{jt}\beta_{j}+v_{jt}+u_{jt},
\end{equation*}
where $C_{jt}$ is the cost of the $j^{\rm{th}}$ bank in year $t$, and $X_{jt}$ contains a constant, the logarithms of the outputs (installment loans, installment loans (to individuals for personal/household expenses), real estate loans, business loans, federal funds sold and securities purchased under agreements to resell, other assets (assets that cannot be properly included in any other asset items in the balance sheet)) and the logarithms of the relevant prices of the inputs (labour, capital, purchased funds, interest-bearing deposits in total transaction accounts, interest-bearing deposits in total non-transaction accounts); relevant prices calculated by dividing total expenses of each input  by the corresponding quantity. The results for the 128 banks are presented in \autoref{datapower3}. The p--values are clearly non--significant, and thus the null hypothesis ${\cal{H}}_0$ of a normal/gamma SFM can not be rejected for the years 1990 and 2000. Nevertheless, the estimated value of the gamma shape parameter is borderline to $p=0$ for the year 1990, and that observation combined with the estimate of the stable tail--index $\widehat \alpha\approx 1.89$ reported in Table ST.2 of the on--line Supplement raises doubts whether a normal/gamma is a good model or if   a sub--Gaussian stable/gamma SFM is more appropriate for this case.

\begin{table*}[h]
\caption{Parameter estimates and p-value of the Normal/Gamma test for the dataset in Tsionas (2006) with COLS estimation}
\label{datapower3}
\centering
\setlength{\tabcolsep}{10pt} 
\renewcommand{\arraystretch}{0.75}
\small
{%
\begin{tabular}{l cccc  l cccc} 
\toprule
\multirow{10}{*}{1990} & $\widehat{\beta}_0$ & 3.4383 & & & \multirow{10}{*}{2000} & $\widehat{\beta}_0$ & 1.0220 &  &   \\
& $\widehat{\beta}_1$ & -0.0419 & $\widehat{\sigma}_v$ & 0.2093 & & $\widehat{\beta}_1$ & 0.4349 & $\widehat{\sigma}_v$ & 0.2163  \\
& $\widehat{\beta}_2$ & 0.0376 & $\widehat{p}$ & 1.8E-07 & & $\widehat{\beta}_2$ & -0.0476 & $\widehat{p}$ & 0.2210 \\
& $\widehat{\beta}_3$ & -0.0168 & $\widehat{c}$ & 5.9827 & & $\widehat{\beta}_3$ & 0.3204 & $\widehat{c}$ & 0.2978 \\
& $\widehat{\beta}_4$ & 0.3089 & \multicolumn{2}{c}{{\bf{p--value}}=1.0} &&  $\widehat{\beta}_4$ & 0.0640 & \multicolumn{2}{c}{{\bf{p--value}}=0.98} \\
& $\widehat{\beta}_5$ & 0.0243 & & & & $\widehat{\beta}_5$ & -0.0093 & &\\
& $\widehat{\beta}_6$ & 0.3003 & & & & $\widehat{\beta}_6$ & 0.3737 & & \\
& $\widehat{\beta}_7$ & 0.1295 & & & & $\widehat{\beta}_7$ & 0.1282 & &\\
& $\widehat{\beta}_8$ & 0.0556 & & & & $\widehat{\beta}_8$ & 0.030 & & \\
& $\widehat{\beta}_9$ & 0.4628 & & & & $\widehat{\beta}_9$ & 0.3642 & & \\
\bottomrule
\end{tabular}}
\end{table*}

The conclusions drawn by the three last applications are in partial agreement with those of  Horrace and Wang \citeyearpar{HW22} obtained by an entirely non--parametric  method, with some discrepancies potentially  be attributed to the lack of fit of the particular law under test. In this connection the stable law advocated herein as a heavy--tailed law is one option while the Student--t distribution may also serve as a potential heavy tailed model for  the pure random error component $v$; see Wheat et al. \citeyearpar{WSG}. The two families however are  generally separate families of distributions  having only the Cauchy law in common, a law that is extremely heavy--tailed  which makes it an unlikely candidate for the distribution of $v$.

We close this section by noting that  with our methods we may be able to distinguish between  a medium--tailed normal/gamma and a purely heavy--tailed stable/gamma SFM.  The practitioner can apply the stable/gamma test in \eqref{test4} as well as the MGF--based test in \eqref{statistic}, and lean towards either model depending on the results, i.e. on the corresponding p--values, and the parameter estimates obtained, in a manner similar to that presented in the real--data examples of this section. Of course judgement should be exercised,  and in any case  the procedure is heuristic in that there exist obvious dependencies of test results applied on the same dataset as well as uncertainty due to variability of parameter--estimates. On a different note and using our MLE,  the identification of a heavy--tailed SFM may also be assessed by a combination of the CF--based test followed by a likelihood--ratio test.  A CF--based test to decide whether (any) stable/gamma SFM is data--compatible, which,  in the case of acceptance, should be followed by a likelihood--ratio test  to see if the restriction $\alpha=2$ is supported by the data.
In this connection we point out  that applying a likelihood--ratio test alone leads to a directional procedure, i.e. it restricts attention of the pool of potential SFM to the stable/gamma class, and in doing so ignores other non--nested models. Therefore successive application of CF--based and likelihood--ratio tests is a more appropriate global procedure.

\section{Conclusions}\label{sec_7}
We propose goodness--of--fit tests in SFMs based on empirical versions of distributional transforms (the moment generating and the characteristic function) of the law of the composed error term $\varepsilon=v-u$. Our tests are somewhat complicated in functional form, and in the case of the stable/gamma test make use of special functions, a price that often has to be paid for an omnibus test with such complicated models. At the same time though the test statistics are free of numerical integration and thus easier to implement compared to classical tests based on the distribution function. Also conditionally on correct specification of the production function (or cost function),  the interpretation of test results depends on the postulated pure random error term: If we are confident enough about the law of $v$ then the test result has a direct and definite bearing on the law of inefficiency component $u$, while otherwise rejection of the null hypothesis may also be due to poor fit of the postulated law of $v$. We study  the finite--sample properties of the tests in a series of Monte Carlo simulations by means of  bootstrap versions of the test statistics. Moreover the application of our methods  on real data  serves the purpose of investigating the presence of heavy--tails in the empirical modelling of production.

\newpage

\appendix
\appendixpage{}

\renewcommand{\thesection}{A.\arabic{section}}

\section{Consistency of the MGF--based test} \label{ap1}
Assume that the law of the composed error $\varepsilon$ has a  finite and differentiable MGF. 
Assume also that $(\widehat \sigma^2_v,\widehat c,\widehat p)\to (\sigma_v^2,c,p) \in (0,\infty)^3$, almost surely,  i.e. that the parameter estimates attain positive and finite almost sure limits even under alternatives. Then we have by the pointwise consistency of the empirical MGF (see Feuerverger, \citeyearpar{Feu89}) that,  $D^2_n(t)\to D^2(t)$, where  $D(t)=((1+t)/c)M'_\varepsilon(t/c)+[p-\lambda t(1+t)]M_\varepsilon(t/c)$.  If we further assume  that the weight function $w>0$ is such that $$\int_{0}^\infty t^4 e^{tx} w(t) {\rm{d}}t<\infty, \ \forall x\in\mathbb R,$$ it then follows by Fatou's lemma (see Jiang \citeyearpar{Jiang}, \S A.2.3) and \eqref{ts} that \[ \liminf_{n\to\infty} \frac{T_{n,w}}{n} \geq  \int_{0}^\infty D^2(t) w(t) {\rm{d}}t:=\Delta_w, \]
whereby a simple change of variables
\[\Delta_w=c^{-1} \int_{0}^\infty \widetilde D^2(t) w(t) {\rm{d}}t,\] 
with 
\[
\widetilde D(t)=(1+c t)M'_\varepsilon(t)+[c p-\sigma^2_{v}t(1+c t)]M_\varepsilon(t).
\]
Clearly $\Delta_w>0$ unless the integrand above vanishes identically in $t$. However since $w>0$ and due to the uniqueness of solution of the differential equation in \eqref{mgf14} (with initial condition $M_\varepsilon(0)=1$),  this integrand vanishes only  if the MGF of $\varepsilon$ is given by  \eqref{mgf15}. In turn, by the uniqueness of MGFs  the last implication is equivalent to the null hypothesis ${\cal{H}}_0$ figuring in \eqref{null}, i.e. to a normal/gamma SFM with the indicated parameters. Otherwise under alternatives,  $T_{n,w}\to \infty$, a.s. as $n\to\infty$, and thus the test which rejects  ${\cal{H}}_0$ for large values of $T_{n,w}$ is consistent.   

\begin{remark}
Consistency of the MGF--based test critically depends on the one--to--one correspondence between distributions and MGFs. However strictly speaking this uniqueness assumes that the  underlying MGF is defined in a neighborhood of zero (see Severini, \citeyearpar{Sev05}, Theorem 4.9). Consequently one needs to integrate in the test statistic \eqref{ts} over an interval containing the origin. Here for convenience we integrate over $(0,\infty)$. A corresponding statistic of the type $n \int_{t_0}^\infty D^2_n(t) e^{-\gamma t^2}{\rm{d}}t$, with $t_0<0$, was also tried in simulations, but the results  produced are very similar and therefore not reported.
\end{remark}

\section{Computation of the MGF--based statistic} \label{ap2}

Starting from equation \eqref{mgf13} we compute
\begin{eqnarray*}
D^2_n(t)&=&(1+t)^2 M'^2_n(t)+[\widehat p-\widehat \lambda t(1+ t)]^2 M^2_n(t) +2(1+t)[\widehat p-\widehat \lambda t(1+ t)] M_n(t)M'_n(t) \\ &=&(1+t)^2 \frac{1}{n^2} \sum_{j,k=1}^n r_j r_k e^{t(r_j+r_k)} +[\widehat p-\widehat \lambda t(1+ t)]^2 \frac{1}{n^2} \sum_{j,k=1}^n e^{t(r_j+r_k)} \\ &+&2(1+t) [\widehat p-\widehat \lambda t(1+ t)]\frac{1}{n^2} \sum_{j,k=1}^n r_j e^{t(r_j+r_k)}.
\end{eqnarray*}
Now plugging   the above expression for $D_n^2(t)$ into the test statistic \eqref{ts} and integrating term-by-term we arrive, after some grouping,  at the reported expression of $T_{n,\gamma}$, by  making use of the integrals 
\[\int_{0}^\infty t^k e^{t x}  e^{-\gamma t^2}dt, \ k=0,1,2,3,4.\]
These integrals may be obtained in closed form by use of the standard normal distribution function $\rm\Phi(\cdot)$. By way of example for $k=0, 1, 2$ we have
\[\int_{0}^\infty  e^{t x}  e^{-\gamma t^2}dt=\sqrt{\frac{\pi}{\gamma}} {\rm{\Phi}}\left(\frac{x}{\sqrt{2\gamma}}\right) e^{\frac{x^2}{4\gamma}},\]
\[\int_{0}^\infty t \:e^{t x}  e^{-\gamma t^2}dt=\frac{1}{2\gamma} +\sqrt{\frac{\pi}{4\gamma^3}} {\rm{\Phi}}\left(\frac{x}{\sqrt{2\gamma}}\right) x \:e^{\frac{x^2}{4\gamma}},\]
and
\[\int_{0}^\infty t^2 \:e^{t x}  e^{-\gamma t^2}dt=\frac{x}{4\gamma^2} +\sqrt{\frac{\pi}{16\gamma^5}} {\rm{\Phi}}\left(\frac{x}{\sqrt{2\gamma}}\right) (x^2+2\gamma) \:e^{\frac{x^2}{4\gamma}},\]
respectively. For $k=3$ and $k=4$ the corresponding formulae are slightly  more complicated but still manageable.

\section{COLS and ML estimation} \label{ap3}

\subsection{COLS estimation and a limit statistic}  
As already mentioned we use the COLS estimator of Kopp and Mullahy \citeyearpar{KM93}. Specifically  estimates of the  regression coefficients are obtained by using the first moment equation in \eqref{moments_stand}, while the estimates of the distributional parameters $(\sigma^2_v,p,c)$ are obtained by using, in addition to the first moment equation ${\cal{M}}_1$, the next three equations, i.e. ${\cal{M}}_2,{\cal{M}}_3$ and ${\cal{M}}_4$, figuring in \eqref{moments_stand} and thereby solving the system of four equations so produced. 

In this connection, and in order to obtain the COLS estimators, we consider the SFM with a corrected error term 
\begin{eqnarray}
Y_j=X_{j}+v-(u-\mathbb E(u)) = X_{j}+v -\delta \text{ , where } \mathbb E(u)=cp , \nonumber
\end{eqnarray} 
and we solve the following system of equations:

\begin{eqnarray}.
& & \frac{1}{n}\sum^{n}_{j=1} X^\top_{j}(Y_j-X_{j}\beta +cp)=0 , \nonumber \\
& & \frac{1}{n}\sum^{n}_{j=1} (Y_j-X_{j}\beta +cp)^2 - (\sigma_v^2 + c^2p)=0 ,\nonumber \\
& & \frac{1}{n}\sum^{n}_{j=1} (Y_j-X_{j}\beta +cp)^3 + 2c^3p=0 , \nonumber \\
& & \text{and} \nonumber \\
& & \frac{1}{n}\sum^{n}_{j=1} (Y_j-X_{j}\beta +cp)^4 - (3\sigma_v^4 + 6 \sigma_v^2 c^2 p +3p(2+p)c^4)=0 . \nonumber
\end{eqnarray}

Referring back to equation \eqref{Dn}, notice that the use of COLS estimators of the parameters, leads to ${\cal{M}}_{k,n}=0, \ k=1,2,3,4$. Therefore inserting  \eqref{Dn} in the test statistic equation \eqref{ts},  and integrating $D^2_n(t)$ term--by--term we deduce
\begin{eqnarray} \nonumber
T_{n,\gamma}&=&n\left(\frac{{\cal{M}}_{5,n}}{24}\right)^2 {\cal{I}}_{8,\gamma}+{\rm{o}}(\gamma^{-9/2}), \ \gamma \to \infty, 
\end{eqnarray}
where 
\[ {\cal{I}}_{k,\gamma}:=\int_0^\infty t^k e^{-\gamma t^2} {\rm{d}}t=\frac{{\rm{\Gamma}}(\frac{k+1}{2})}{2\gamma^{\frac{k+1}{2}}},\] and consequently \begin{eqnarray} \nonumber
\lim_{\gamma \to \infty} \frac{2\gamma^{9/2}}{{\rm{\Gamma}}(9/2)}  T_{n,\gamma}&=&n \left(\frac{{\cal{M}}_{5,n}}{24}\right)^2.  
\end{eqnarray}
The last equation  shows that conditionally on COLS estimation, the test statistic  $T_{n,\gamma}$ is dominated by the fifth empirical moment equation,  in the limiting case $\gamma \to \infty$.  This might appear a bit confusing in view of the remarks made at the end of Section \ref{sec_2} that higher values of $\gamma$ favor lower order moments, but it should be kept in mind that in this case the moment of order five is the lowest non-vanishing moment. We finally note that, on account of  equation \eqref{moments_stand} and as $n\to\infty$, this empirical moment equation vanishes under the null hypothesis ${\cal{H}}_0$ of a SFM with a normal/gamma composed error distribution.

\subsection{Likelihood Computation} 

The test procedure based on $T_{n,\gamma}$, and its counterpart 
for the stable/gamma model (to be presented in the next 
section), does not necessarily depend on COLS estimation since the estimated parameters may be obtained by maximum likelihood with an extra computing cost. The extra cost is due to the fact that the inefficiency term $u$ 
is a latent variable which must be integrated out when
computing the likelihood. Roughly speaking, the idea is to maximize the likelihood, whereby the density is obtained  by  Fourier inversion from the corresponding CF. Except in a few instances, this process  
will require numerical integration as is the case in both the normal/gamma
and stable/gamma models. Following 
Tsionas \citeyearpar{Tsionas2012} and Andrade and Souza \citeyearpar{Andrade2018}
we approximate the density of the combined error, $f_\varepsilon$, 
by numerical inversion of its CF, 
\begin{equation}
\label{eq:inv}
f_\varepsilon(\varepsilon) 
= 
\frac{1}{2\pi} \int_{-\infty}^{\infty}
\varphi_\varepsilon(t) \exp(-{\rm{i}}t\varepsilon) {\rm{d}}t.
\end{equation} 
The CF in the normal/gamma model is given by
\begin{equation}\label{eq:cfNG}
\varphi_\varepsilon(t) = \exp \left( 
\frac{-t^2\sigma_v^2}{2}
\right)
(1+{\rm{i}}ct)^{-p}, \ t \in \mathbb{R}, 
\end{equation}
and, in the stable/gamma case, by 
\begin{equation}\label{eq:cfSG}
\varphi_\varepsilon(t) = \exp \left( 
-\kappa^\alpha |t|^\alpha
\right)
(1+{\rm{i}}ct)^{-p}, \ t \in \mathbb{R}.
\end{equation}

The approximated density at $N$ (typically large) evenly spaced 
values $e_k$ and $t_j$ is 
\begin{equation}\label{eq:dft}
\hat{f}_\varepsilon(e_k)  
= a_N \mathrm{IDFT}(z_j) ,
 \end{equation}
where $a_N$ is a term depending on $N$ and the spacing used, and
IDFT is the (inverse) discrete Fourier transform of the 
sequence $z_j = (-1)^j \varphi_\varepsilon(t_j)$ 
which can be efficiently calculated by the
fast Fourier transform algorithm. Details of the approximation
of \eqref{eq:inv} by \eqref{eq:dft} are given in 
Andrade and Souza \citeyearpar{Andrade2018} with a very similar notation.

\subsection{Efficiency Measures} \label{sec:eff-meas}

Another direct use of the likelihood approximation above is to calculate 
efficiency measures. Given parameter estimates, 
we can obtain the composed 
errors $\hat{\varepsilon} = y - \hat{y}$. 
With a single estimate $\hat{\varepsilon}_i$ for the $i$-th firm, we cannot estimate $u_i$ and $\nu_i$ separately. However,
this is not necessary for the main 
purpose of efficiency analysis which is
the estimation of individual
efficiency scores. The main ones are
$\mathbb E(\exp(-u)|\varepsilon)$ by Battese and Coeli \citeyearpar{BC1995} and $\exp(-\mathbb E(u|\varepsilon))$ following Jondrow et al.~\citeyearpar{JLMS} 
and, less often, Mode$(u|\varepsilon)$.

There are no explicit formulas for these measures except for 
the simplest SFMs such as normal/exponential or 
normal/half-normal but the 
techniques just presented for the
computation of the likelihood can be use in the normal/gamma 
and stable/gamma models. 

For instance, 
$$
\mathbb E(e^{-u}|\varepsilon=z) \approx 
\frac{1}{\hat{f}_\varepsilon(z)}
\int_0^{\infty}
\exp(-u)f_u(u) \hat{f}_\nu(z+u) \mathrm{d}u.
$$
Note that $\hat{f}_\varepsilon$ and
$\hat{f}_\nu$ are denoted with a 
hat, indicating that they are approximated by 
the inversion of the characteristic function proposed above as in expression \eqref{eq:dft}.  In the normal/gamma case, 
$f_\nu$ can be obtained explicitly but $f_\varepsilon$ still needs
to be approximated.

\section{Computation of the CF--based statistic} \label{ap4}
We follow analogous steps as those in \ref{ap2}. Specifically we plug  the expression for $\Delta_n^2(t)$ given by  
\eqref{d2} into the test statistic \eqref{test3} along with the formulae
\[ |\varphi_n(t)|^2=\frac{1}{n^2} \sum_{j,k=1}^n \cos t(r_j-r_k) , \]
\[ |\varphi'_n(t)|^2=\frac{1}{n^2} \sum_{j,k=1}^n r_j r_k \cos t(r_j-r_k) , \]
\[ S'_n(t) C_n(t)-C'_n(t)S_n(t)=\frac{1}{n^2} \sum_{j,k=1}^n r_j \cos t(r_j-r_k) , \] 
\[ C'_n(t) C_n(t)+S'_n(t)S_n(t)=\frac{1}{n^2} \sum_{j,k=1}^n r_k \sin t(r_j-r_k) , \] 
and then integrate term-by-term to get to \eqref{test4}.  

\section{Integrals involved in the CF--test} \label{ap5}
The  integral $I_{\nu,\gamma}$ may be computed as
\[
I_{\nu,\gamma}(z)=\frac{1}{2\gamma^{\frac{\nu+1}{2}}}{\rm{\Gamma}}\left(\frac{\nu+1}{2}\right) \ e^{-z^2/4\gamma} \  _{1}F_{1}\left(-\frac{\nu}{2};\frac{1}{2};\frac{z^2}{4\gamma}\right),
\]
while the integral $J_{\nu,\gamma}$ as
\[
J_{\nu,\gamma}(z)=\frac{z}{2\gamma^{1+\frac{\nu}{2}}}{\rm{\Gamma}}\left(1+\frac{\nu}{2}\right) \ e^{-z^2/4\gamma} \ _{1}F_{1}\left(\frac{1-\nu}{2},\frac{3}{2},\frac{z^2}{4\gamma}\right),
\]
where $_1$$F_{1}$ stands for the Kummer confluent hypergeometric function; see Gradshteyn and Ryzhik \citeyearpar{GR94}, 3.952, eqn. 7 and 8 (p. 529), and \S9.2 (p. 1084).

\end{document}